\documentclass[12pt,letterpaper]{amsart}
\usepackage{geometry}
\usepackage{mathrsfs,amsmath,amssymb}
\usepackage{amsfonts,mathpazo}
\usepackage{amsthm}
\usepackage{setspace}
\newtheorem{theorem}{Theorem}
\newtheorem{lemma}{Lemma}
\newtheorem{corollary}{Corollary}

\theoremstyle{remark}

\newtheorem{definition}{Definition}
\newtheorem{example}{Example}

\newcommand{\gd}{\Gamma_{\D}}
\newcommand{\gdtwo}{\Gamma_{\D_{2}}}
\newcommand{\gdone}{\Gamma_{\D_{1}}}
\newcommand{\C}{\mathbb{C}}
\newcommand{\Lev}{\mathcal{L}}
\newcommand{\D}{\Omega}
\newcommand{\Wb}{\overline{W}}

\newcommand{\Dc}{\overline{\Omega}}
\newcommand{\dbar}{\overline{\partial}}
\newcommand{\p}{(P)}
\newcommand{\ptilde}{(\widetilde{P})}

\begin{document}
\onehalfspace
\title{Strong Stein Neighborhood Bases}
\author{S\"{o}nmez \c{S}ahuto\u{g}lu}

\address{University of Toledo, Department of Mathematics, Toledo, OH 43606, USA}
\email{Sonmez.Sahutoglu@utoledo.edu}
\thanks{This is a revision of the preprint formerly titled ``Holomorphic Invariance of Stein
Neighborhood Basis''}
\thanks{Research supported in part by NSF grant numbers DMS-0602191 and DMS-0500842,  and  the Erwin
Schr\"{o}dinger International Institute for Mathematical Physics.}
\thanks{ 2000 \emph{Mathematics Subject Classification}: Primary 32W05; Secondary 32A38}
\keywords{Stein neighborhood basis, $\overline{\partial}$-Neumann problem, pseudoconvex domains}
\date{\today}
\begin{abstract}
Let $\D$ be a smooth bounded pseudoconvex domain in $\C^n .$ We give several 
characterizations for the closure of $\D$ to have  a strong Stein neighborhood
basis in the sense that $\D$ has a defining function $\rho$ such that 
$\{z\in \C^{n}:\rho(z)<\varepsilon\}$ is pseudoconvex for sufficiently
small $\varepsilon>0$. We also show that this condition is invariant under
proper holomorphic maps that extend smoothly up to the boundary.
\end{abstract}

\maketitle

\section{introduction}

A domain  $\D\subset\C^{n}$ is called a domain of holomorphy if there exists a holomorphic function
on $\D$ that cannot be ``extended'' past any boundary point. Any domain in $\C$ is a domain of
holomorphy. However, Hartogs in 1906 discovered that  not every domain in $\C^{n}$ for $ n\geq 2$ is
a domain of holomorphy. This fundamental discovery led to the characterization of domains of
holomorphy called the Levi problem. The Levi problem was first solved by Oka in 1930's for $n=2$ and
by Bremermann, Norguet, and Oka in 1950's for $n\geq 3$. The solution of the Levi problem revealed a
very interesting fact about  domains of holomorphy: domains of holomorphy are precisely the
so-called pseudoconvex domains and hence can be exhausted by  pseudoconvex subdomains. That is, one
can ``approximate" a  domain of holomorphy (or a pseudoconvex domain) from inside by pseudoconvex
domains. Therefore, it is natural to ask whether it is possible to approximate such domains from
outside. We refer the reader to \cite{Ho,Kr,Ra} for precise definitions and basic facts about
domains of holomorphy and pseudoconvex domains.

A compact set $K\Subset \C^{n}$ is said to have a Stein neighborhood basis if for any domain   $V$
containing $K$ there exists a pseudoconvex domain $\D_V$ such that $K\subset \D_V \subset V$. It is
worth noting that  the closure of the Hartogs triangle, $\{(z,w)\in \C^2:0\leq |z|<|w|<1\}$, does
not have a Stein neighborhood basis. However, the Hartogs triangle is not smooth. In 1977, 
Diederich and Forn\ae{}ss(\cite{DF77a}) found a smooth bounded pseudoconvex domain, the so-called
worm domain, whose closure does not have a Stein neighborhood basis, thus answering in the negative
a question of Behnke and Thullen(\cite{BehThullen}). Subsequently, the question of when a Stein
neighborhood basis exists has been studied by Bedford and Forn{\ae}ss  \cite{BF78}, Diederich and
Forn{\ae}ss \cite{DF77a,DF77b}, Sibony \cite{Sib87,Sib91}, Stens{\o}nes \cite{Sten}, et al. 

The existence of a special kind of a  Stein neighborhood basis is known to be connected to global
regularity of the  $\overline{\partial}$-Neumann problem \cite{St01}, approximation properties for
holomorphic functions  \cite{Cir,FN}, and uniform algebras \cite{Ros61}. 

In this paper we will concentrate on smooth domains and  ``smooth'' means 
$C^{\infty}$-smooth. However, the reader will notice that some of the results are still true for
domains with $C^{3}$-smooth boundary. We are interested in the following stronger notion of Stein
neighborhood bases for smooth domains, as it is fairly general and has many applications. We say the
closure $\Dc$ of a smooth bounded pseudoconvex domain $\D$ has a {\it{strong Stein neighborhood
basis}} if $\D$ has a defining function $\rho$ (see the Section \ref{sec2} for a definition) and
there exists $\varepsilon_{0}>0$ such that $\{z\in \C^{n}: \rho(z)<\varepsilon\}$ is pseudoconvex
for $0\leq \varepsilon\leq \varepsilon_{0}$. In Theorem \ref{MainThm} we give several
characterizations for the closure to have a strong Stein neighborhood basis. The precise statement
of Theroem \ref{MainThm} requires some technical definitions, and so is postponed to Section
\ref{sec2-1}. We note that all smooth bounded pseudoconvex domains whose closure known to have a
Stein neighborhood basis satisfy this condition. Whether it is equivalent to having a Stein
neighborhood basis for the closure of a smooth bounded pseudoconvex domain is still open. 

The existence of strong Stein neighborhood bases  implies the so-called uniform H-convexity. A
compact set $K\subset \C^{n}$ is said to be uniformly $H$-convex if there exists a positive sequence
$\{\varepsilon_{j}\}$ that converges to 0, $c>1$, and a sequence of pseudoconvex domains $\D_{j}$
such that $K\subset \D_{j}$ and $\varepsilon_{j}\leq {\rm dist}(K,\C^{n}\setminus\D_{j})\leq
c\varepsilon_{j}$ for $j=1,2,\ldots$. {\v{C}}irka(\cite{Cir}) showed that uniform H-convexity
implies a ``Mergelyan-like'' approximation property for holomorphic functions. There are three
conditions that are known to imply the existence of a (strong) Stein neighborhood basis for the
closure of a smooth bounded pseudoconvex domain in $\C^{n}:$ having a holomorphic vector field in a
neighborhood of the weakly pseudoconvex points  that is transversal to the boundary \cite{FN},
property $\ptilde$ \cite{Sib87,Sib91}, and having a defining function that is
plurisubharmonic on the boundary \cite{FH,FH2}.

The following example shows that having a Stein neighborhood basis for the closure of a domain is
not a invariant under biholomorphism in general. 

\begin{example}
Let $\D_{1}=\{(z,w)\in \C^{2}:0\leq |z|<|w|<1\}$ be the Hartogs triangle, and
$\D_{2}=\{(z,w)\in \C^{2}:0\leq|z|<1,0<|w|<1\}.$ Let $F:\D_{1}\to \D_{2}$ be a
biholomorphism defined as follows: $F(z,w)=(z/w,w).$ One can check that
$\overline{F(\D_{1})}=\Dc_{2}.$ Therefore, although $\overline{F(\D_{1})}$ has a
Stein neighborhood basis, $\Dc_{1}$ does not. 
\end{example}

It is still open whether having a Stein neighborhood basis for the closure of a smooth bounded
pseudoconvex domain is invariant under biholomorphisms in general. However, an easy corollary to our
main result (Theorem \ref{MainThm}), is that having a strong Stein neighborhood basis is invariant
under biholomorphisms that extend up to the boundary.  More precisely, 

\begin{corollary}\label{Invariance}
Let $\D_1$ and $\D_2$  be two smooth bounded pseudoconvex domains in $\C^n, n\geq 2$. Assume that
there exist a proper holomorphic map $F:\D_1\to \D_2$ that extend smoothly to $\Dc_1$ and $\Dc_{2}$
has a strong Stein neighborhood basis. Then  $\Dc_1$ has a strong Stein neighborhood basis.
\end{corollary}

At this point we would like to mention two open questions:  Is the assumption of smooth
extendibility of the biholomorphism in the Corollary above needed, or is it automatic? Does having a
Stein neighborhood basis for the closure of a smooth bounded pseudoconvex domain imply the existence
of a strong Stein neighborhood basis for the closure of the domain?
 
This paper is organized as follows: in Section \ref{sec2} we  set the notation and give basic
definitions.  In Section \ref{sec2-1} we state  the main theorem, Theorem \ref{MainThm}, and give an
application to a potential theoretic property the so-called property $\ptilde$ (see Corollary
\ref{P-tilde}). In Section \ref{sec3} we give the proof of Corollary \ref{Invariance} and Theorem
\ref{MainThm}. In Section \ref{sec4} we give the proof of Corollary \ref{P-tilde}.

 We would like to take this opportunity to thank David Barrett, Harold Boas, Mehmet \c{C}elik, John
Erik Forn\ae{}ss, Daniel Jupiter, Berit Stens{\o}nes, and my advisor Emil Straube for reading the
early manuscripts, valuable comments,  and stimulating discussions. This article is based on a part
of the author's Ph.D. thesis \cite{Sah06}.
 
\section{Notation and Definitions}\label{sec2}

Let $\D$ be a smooth bounded pseudoconvex domain in $\C^{n}$ and $r$ be a defining function for
$\D$. That is, $r$ is a smooth function defined in a neighborhood $U$ of $\Dc$ such that it is
negative on $\D$, positive on $U\setminus \Dc$, zero on the boundary $b\D$ of $\D$, and the gradient
$\nabla r$ of $r$ does not vanish on $b\D$. We define the complex Hessian of $r$ at 
$z$ as follows: 
$$\Lev_r (z;A,\overline{B})=\sum_{j,k=1}^n\frac{\partial^2 r(z)}{\partial z_j\partial \bar
z_k}a_j\bar b_k,$$  
where $A=(a_1,\ldots,a_n)$ and $B=(b_1,\ldots,b_n)$ are vectors in $\C^n .$  We would like to note
that  we identify $\C^n$ with the $(1,0)$ tangent bundle of $\C^n$. Namely, $(a_1,\ldots,a_n)$ is
identified with $\sum_{j=1}^n a_j\frac{\partial}{\partial z_j}.$ 
We will denote $\Lev_r (z;A,\overline{A})$ by $ \Lev_r (z;A),$ and   $\sum_{j=1}^n|w_j|^2$ by
$\|W\|^2,$ where $W=(w_1,\ldots ,w_n)\in \C^n .$  Let $\vec{n}(z)\in \mathbb{R}^{2n}$ be the unit
outward normal vector of $b\D$ at $z.$ We denote the  directional derivative in the direction
$\vec{n}(z)$ at the point $z$ by $\frac{d}{d \vec{n}(z)}.$  To simplify the notation and ease the
calculation we will use the following notation: 
$$\frac{d}{dr(z)} =\frac{1}{\|\nabla r(z)\|} \frac{d}{d\vec{n}(z)}\,\,\, {\rm and}\,\,\,
A(h)(z)=\sum_{j=1}^n \frac{\partial h (z)}{\partial z_j}a_j
$$
for a (type $(1,0)$) vector $A=(a_1,\ldots,a_n)\in \C^n $ and $h\in C^1(\overline{\D}).$ 
It is a standard fact that a smooth bounded domain $\D$ is pseudoconvex if and only if
$\Lev_{r}(z;W)\geq 0$ for $z\in b\D$ and $W(r)(z)=0$ where $r$ is a defining function for $\D$. One
can check that this condition is independent of the defining function $r$. When we refer to finite
or infinite type of a point in $b\D$, we mean type in the sense of D'Angelo \cite{D}. 
Let $\D_{\infty}$ denote the set of infinite type points of $b\D$ and
$$ \gd=\{(z,W)\in b\D\times \C^n: z \in \D_{\infty},W(r)(z)=\Lev_{r}(z;W)=0,\|W\|=1\}.$$
$\gd$ is, in some sense, the unit sphere of the weakly pseudoconvex directions on the infinite type
points. 
For a fixed vector $A\in \C^n$ and $z\in b\D$ we will denote 
\begin{eqnarray*}C_r(z;A)&=&\frac{d\Lev_r (t;A )}{dr(z)}\bigg|_{t=z}, \\
                      D_r(z;A)&=&\Lev_r(z;A,\overline N_r),\,\, {\rm and}\\
		      E_{r}(z;A)&=&C_{r}(z;A)-2{\rm Re}\bigg(D_r(z;A)\overline{A}(\ln \| \nabla
r\|)(z)\bigg),
\end{eqnarray*}
where $r$  is a defining function for $\D$ and $N_r=\frac{4}{\| \nabla r(z)\|^2
}\sum_{j=1}^n\frac{\partial r}{\partial \bar z_j}\frac{\partial}{\partial z_j}.$ Notice that
$C_r(z;A)$ is the directional derivative of the complex Hessian $\Lev_r(z;A)$ in the (real) normal
direction at $z$ for a fixed vector $A$. $C_r$ alone does not guarantee that the level sets outside
the domain will be pseudoconvex because it measures how the complex Hessian changes as one moves out
of the domain for only  fixed vectors. On the other hand, $E_{r}(z;A)$ gives a sufficient condition
for the existence of a Stein neighborhood basis (see Theorem \ref{MainThm}) because it takes into
account how the complex tangent vectors change as one moves out of the domain. This change can be
measured by $D_{r}(z;A)$. We   note that $D_r(z;A)$ plays a very crucial role in the vector field
approach of Boas and Straube for the global regularity of the $\dbar$-Neumann problem
\cite{BS91,BS93,BS}. This might suggest that there are deeper relations between the global
regularity of the $\dbar$-Neumann problem and the existence of a Stein neighborhood basis for the
closure. 

\section{Statement of the Main Theorem}\label{sec2-1}

The following is our main theorem. It gives several characterizations of having a strong Stein
neighborhood basis for the closure of a smooth bounded pseudoconvex domain. It will be used in the
proof of Corollary \ref{Invariance}.
  
\begin{theorem}\label{MainThm}
Let $\D$ be a smooth bounded  pseudoconvex domain in $\C^n ,n\geq 2.$ The following conditions are
equivalent:

\begin{itemize}
\item[i)] there exist a neighborhood $U$ of $\Dc$, a defining function $\rho$ for $\D$ in $U,$ and
$c>0$ such that $\Lev_{\rho}(z;W)\geq c\rho(z)\|W\|^2$ for $z\in U\setminus \D$ and $W(\rho)(z)=0,$
 \item[ii)] there exists a  defining function $\rho$ for $\D$ and $\varepsilon_{0}>0$ such that $\{
z\in \C^n:\rho(z)<\varepsilon \}$ is pseudoconvex  for  $0\leq\varepsilon\leq \varepsilon_{0}$. That
is, $\Dc$ has a strong Stein neighborhood basis
\item[iii)] there exists a  defining function $\rho$ for $\D$ such that
$$E_{\rho}(z;W)\geq 0 \,\,for\,\,(z,W)\in \gd .$$
\item[iv)] there exist $h\in C^{\infty}(\overline{\D})$  and a  defining function $r$ for $\D$  such
that:
 \begin{equation}\label{eqP_0}\displaystyle \Lev_h(z;W)\geq |W(h)(z)|^2+2{\rm Re}\Big(
D_r(z;W)\overline{W}(h)(z)\Big)-E_r(z;W) \end{equation}
for $(z,W)\in \gd .$ 
\end{itemize}
 
\end{theorem}

Now we will give the definition of a potential theoretic condition: property $\ptilde$. The
following definition is from \cite{McN}.
\begin{definition}\label{DefPtildeP} 
Let $\D$ be a bounded domain in $\C^{n}$. Then  $\D$ satisfies property $(\widetilde{P})$  if there
exists a sequence of plurisubharmonic functions $\{\phi_{j}\} \subset C^{\infty}(\Dc)$ such that 
\begin{itemize}
\item [i)] $ |W(\phi_{j})(z)|^2\leq \Lev_{\phi_{j}}(z;W)$ for $z\in \Dc$ and $W\in \C^{n}$,
\item [ii)] $\Lev_{\phi_{j}}(z;W)\geq j\|W\|^{2}$ for $z\in b\D $ and $W\in \C^{n}$.
\end{itemize}
\end{definition}

We note that  if $\D$ is a $C^{3}$-smooth bounded pseudoconvex domain that satisfies property
$\ptilde$ then for any defining function $r$ there exists $h\in C^{2}(\overline{\D})$ such that
\eqref{eqP_0} is satisfied \cite{Sah06}. 
Property $\ptilde$ is  related to another potential theoretic property called property $\p$. The
difference is that instead of i) in the above definition,  property $\p$ requires the sequence
$\{\phi_{j}\}$ to be uniformly bounded on $\Dc$. By exponentiating and scaling one can easily show
that property $\p$ implies property $\ptilde$. However, it is still open whether the converse is
true. Although these properties were introduced for studying compactness of the $\dbar$-Neumann
problem (\cite{Cat,McN}) they naturally appear in the study of Stein neighborhood bases.  We note
that Harrington(\cite{Har06}) showed that property $\p$ implies existence of a Stein neighborhood
basis for the closure when the domain is $C^1$-smooth bounded and pseudoconvex.

As a result of our method we get a characterization of property $\ptilde$ in terms of existence of
strong Stein neighborhood bases in some sense.   

\begin{corollary}\label{P-tilde}
Let $\D$ be a smooth bounded  pseudoconvex domain in $\C^n ,n\geq 2.$ Then $\D$ satisfies property
$\ptilde$ if and only if for every $M>0$ there exists a defining function $\rho$ such that
$E_{\rho}(z;W)>M $ for $(z,W)\in \gd$. 
\end{corollary}

An immediate  implication of the above corollary is that property $\ptilde$ can be localized onto
weakly pseudoconvex directions on infinite type points. We note that the localization of property
$\ptilde$ has been obtained by \c{C}elik(\cite{Celik})  in his thesis before.  

\section{Proof of Theorem \ref{MainThm} and Corollary \ref{Invariance}}\label{sec3}

The following Lemmas will be useful in the proof of Theorem \ref{MainThm}.
\begin{lemma}\label{geq0}
Let $\D$ be a smooth bounded  pseudoconvex domain in $\C^n ,n\geq 2.$ Assume that $\D$ has a
defining function $r$ such that $E_{r}(z;W)\geq 0$ for $(z,W)\in \gd$. Then there exists a defining
function $\rho$ for $\D$ such that $E_{\rho}(z;W)>0$ for $(z,W)\in \gd .$ 
\end{lemma}
\begin{proof} Without loss of generality we may assume that $\Dc$ is contained in the ball centered
at the origin and of radius $\tau={\rm diam}(\D)$. 
Let us define $\rho(z)=r(z)e^{h(z)}$ where 
$$h(z)=\frac{e^{\beta(|z|^2-\tau^2)}}{2\beta\tau^2},$$
and  $\beta=4\sup\{1+|\Lev_{r}(z;N_{r},\overline{W})|^2:(z,W)\in \gd \}.$  Now we would like to
calculate $E_{\rho}$ in terms of $r$. We note that $\frac{d}{d\rho}=e^{-h}\frac{d}{dr}$ and
$N_{\rho}=e^{-h}N_{r}$ on $b\D$. One can check that 
\begin{eqnarray*}
D_{\rho}(z;W)&=&D_{r}(z;W)+W(h)(z)\\
\Wb (\ln\|\nabla \rho \|)(z)&=& \Wb (\ln\|\nabla r \|)(z)+\Wb (h)(z)
\end{eqnarray*}
We note that 
\begin{equation}
\label{forC_r}\frac{dW(r)}{dr(z)}(z) = W( \ln \| \nabla r \|)(z)
\end{equation}
because $W$ is a fixed vector and 
\begin{eqnarray*}
\frac{d  W(r) }{dr(p)}(z)\bigg|_{z=p}&=&\frac{d W}{dr(p)}(r)(p)
+\sum_{j=1}^nw_j\frac{\nabla r(p)}{\| \nabla r(p) \|^2}\cdot \nabla \left(
\frac{\partial r}{\partial z_j}\right)(p)  \\
&=&\frac{d W}{dr(p)}(r)(p)+\frac{1}{2}W\left( \ln \| \nabla r \|^2\right)(p)
\end{eqnarray*} 
We note that the second term in the first equality consists of a summation of
dot product of vectors. Using \eqref{forC_r} one can calculate that
$$
C_{\rho}(z;W)=C_{r}(z;W)+\Lev_h(z;W)+|W(h)(z)|^2+2{\rm Re}\bigg(W(h)(z)\Wb(\ln \|\nabla r\|)
(z)\bigg).
$$
If we put the above calculations together we get 
\begin{equation}\label{E_rho}
E_{\rho}(z;W)= \Lev_h(z;W)-|W(h)(z)|^2-2{\rm Re}\bigg(D_r(z;W)\overline{W}(h)(z)\bigg)+E_{r}(z;W).
\end{equation}
So we only need to show that $$g(z,W)=\Lev_h(z;W)-| W(h)(z)|^2 -\sqrt{\beta}|W(h)(z)|>0\,\,{\rm
for}\,\, (z,W)\in \gd .$$ Let us denote $\sum_{j=1}^nw_j\bar z_j$ by $\langle W,z\rangle$. Then one
can calculate that 
\begin{eqnarray*}
|W(h)(z)|&=& \frac{e^{\beta(|z|^2-\tau^2)}}{2\tau^2}\left|\langle W,z\rangle\right| , \,\,{\rm
and}\\
\Lev_h(z;W)&=&\frac{e^{\beta(|z|^2-\tau^2)}}{2\tau^2}\left(\|W\|^2+\beta \left|\langle
W,z\rangle\right|^2\right) .
\end{eqnarray*}
Therefore, if we use the inequality $\sqrt{\beta}\left| \langle W,z\rangle \right|\leq \beta\left|
\langle W,z\rangle \right|^2+1/4$,  the assumption that $\Dc$ is contained in the ball centered at
the origin and of radius $\tau$, and the fact that $\|W\|=1$ we get
\begin{eqnarray*}
g(z,W)&\geq &  \frac{e^{\beta(|z|^2-\tau^2)}}{2\tau^2} \left( 1+\left(\beta
-\frac{1}{2\tau^2}\right) \left|\langle W,z \rangle \right|^2
-\sqrt{\beta}\left| \langle W,z\rangle \right| \right) \\
& \geq&  \frac{e^{\beta(|z|^2-\tau^2)}}{2\tau^2} \left( \frac{3}{4}-\frac{\|z\|^2}{2\tau^2} \right)
.\end{eqnarray*}
Again since $\Dc$ is contained in the ball centered at the origin and it is of radius $\tau$ we have
$$\frac{3}{4}-\frac{\|z\|^2}{2\tau^2}>0\,\,{\rm for}\,\, z\in \Dc .$$ 
This completes the proof of Lemma \ref{geq0}.
\end{proof}
\begin{lemma}
 \label{LocalP}
Let $\D$ be a smooth bounded  pseudoconvex domain in $\C^n ,n\geq 2$ and $K$ be a compact subset of
$b\D .$ Assume that $z$ is of finite type for every $z\in K$  and  $h\in C^{\infty}(\Dc)$ is given.
Then for every $j>0$ there exists $h_j\in C^{\infty}(\Dc)$ such that $|h_j-h|\leq 1/j$ uniformly on
$\Dc$ and $\Lev_{h_j}(z;W)\geq j\|W\|^2$ for $z\in K$ and $W\in \C^n.$
\end{lemma}
\begin{proof}
 Using Proposition 3 in \cite{Sib87a}  we can construct a smooth finite type pseudoconvex subdomain
$D$ such that $K\subset bD\cap b\D .$  Similar construction is used in \cite{Bell86}.
Catlin(\cite{Cat}) showed that finite type domains satisfy property $(P)$. So  $D$ satisfies
property $\p$. That is, there is a sequence of functions $f_j\in C^{\infty}(\overline{D})$ such that
$1/2 \leq f_j\leq 3/4$ on $\overline{D}$ and $\Lev_{f_j}(z;W)\geq j^2\|W\|^2$ for $z\in bD$ and
$W\in \C^n$. Let $\tilde f_j$ denote a smooth extension of  $f_j$'s to $\Dc$ such that $0\leq \tilde
f_j\leq 1$ on $\Dc$. We finish the proof by choosing $h_j=h+\frac{\tilde f_{j+k}}{j}$ for
sufficiently large $k$.
\end{proof}

\begin{proof}[Proof of Theorem \ref{MainThm}] We note that 
i)$\Rightarrow$ ii) is trivial and  using \eqref{E_rho} iii) $\Rightarrow$ iv) is easy to see. To
prove ii) implies iii) let us assume that $\D$ has a defining function $\rho$ and there exists
$\varepsilon_{0}>0$ such that $\{z\in \C^{n}:\rho(z)<\varepsilon\}$ is pseudoconvex for $0\leq
\varepsilon\leq \varepsilon_{0}.$ Now we would like to differentiate $\Lev_{\rho}(z;W)$ in the
outward normal direction for $(z,W)\in \gd$. If we apply $\frac{d}{d\vec{n}(p)}$ to
$\Lev_{\rho}(z;W)$  for any smooth vector field $W$ of type $(1,0)$ such that $W(\rho)=0$ on a
neighborhood of $\D_{\infty}$ and $(p,W(p)) \in  \gd $ (calculations are similar to the ones in the
proof of Lemma \ref{geq0}) we get 
\begin{equation*}
 \frac{d \Lev_{\rho}(z;W)}{d\rho(p)}\bigg|_{z=p}=C_{\rho}(p;W(p))-2{\rm
Re}\bigg(D_{\rho}(p;W)\overline{W}(\ln \| \nabla \rho\|)(p)\bigg)=E_{\rho}(p;W(p)).\end{equation*} 
 So ii) implies that the right hand side of the above equality is nonnegative. Therefore,  
 $E_{\rho}(p;W(p))\geq 0$ for $(p,W(p))\in \gd .$  

Let us prove that iv)$\Rightarrow$ i): We divide the proof into two parts. In the first part, we
will produce a defining function whose sublevel sets are pseudoconvex (from $\D$'s side) outside of
$\D$ in a neighborhood of the set of infinite type points. In the second part, using property $(P)$,
we will modify this defining function away from infinite type points to get  a strong Stein
neighborhood basis for the closure.

\noindent \textit{Analysis on Infinite Type Points.} 
Using Lemma \ref{geq0} we can  assume that $\D$ has
a defining function $r$ and there exists a function $h\in C^{\infty}(\overline{\D})$
 such that 
$$\Lev_h(z;W)> |W(h)(z)|^2+2{\rm Re}\bigg( D_r(z;W)\overline{W}(h)(z)\bigg)-E_r(z;W)$$ 
for $(z,W)\in \gd .$
We extend $h$ to $\C^n$ as a smooth function and call the extension
$h$. We scale $h$, if necessary, so that there is a neighborhood of $\Dc$ on
which the conditions of the theorem are still satisfied. We define $\rho(z)=
r(z)e^{ h(z)} .$ We will show that there exists a neighborhood $V$ of $\D_{\infty}$, the set of
infinite type points in $b\D$, such that $\nabla \rho $ is nonvanishing on $V,$ and the complex
Hessian of $\rho$ is nonnegative on vectors complex tangential to the level sets of $\rho$ in
$V\setminus \D .$ Since $\D$ is
bounded and $\| \nabla \rho \|$ is continuous and strictly positive on $b\D ,$
the first part of the above argument follows immediately. It suffices to argue near a boundary point
$q$ because $\D$ is bounded.

Let $z\in \C^{n}\setminus \D$, and $ W\in \C^n ,$ be a complex
tangential vector to the level set of $\rho$ at $z$. Namely,
\begin{equation}W(\rho )(z)= e^{h(z)}\bigg( W(r)(z)+
r(z)W(h)(z)\bigg) =0.\label{eqnlevistep2} \end{equation} 
Now we will calculate the complex Hessian of $\rho$ at $z$ in the direction $W$.  So
first we differentiate $\rho$ with respect to $z_j$ to get
\begin{equation}\label{eqnlevistep}\frac{\partial \rho}{\partial z_j}(z)=
e^{h(z)}\frac{\partial r}{\partial z_j}(z)+  r(z) 
e^{h(z)}\frac{\partial h}{\partial z_j}(z)\end{equation}
and if we differentiate (\ref{eqnlevistep}) with respect to $\bar z_k$ we get
\begin{eqnarray*}
\frac{\partial^2 \rho}{\partial \bar z_k\partial z_j}(z) &=&
e^{h(z)}\frac{\partial^2r}{\partial \bar z_k\partial z_j}(z)
+ e^{ h(z)}\frac{\partial h}{\partial \bar z_k}(z)
\frac{\partial r}{\partial z_j}(z)  
+ e^{ h(z)}\frac{\partial r}{\partial \bar z_k}(z)
\frac{\partial h}{\partial z_j}(z) \\
&+ &r(z) e^{ h(z)}
\frac{\partial h}{\partial \bar z_k}(z)\frac{\partial h}{\partial z_j}(z) 
+  r(z) e^{ h(z)}
\frac{\partial^2 h}{\partial \bar z_k \partial z_j}(z) .
\end{eqnarray*}
Using (\ref{eqnlevistep2}) in the last equality we get  
$$   \Lev_{\rho}(z;W) = e^{ h(z)}\bigg(  r(z)
\Lev_h(z;W)+\Lev_r(z;W)-r(z)| W(h)(z)|^2 \bigg) .$$ 
We would like to show that there exists a neighborhood $V$ of infinite type points $\D_{\infty}$
such that 
\begin{equation}\label{eqf}f(z,W)= r(z) \Lev_h(z;W)+\Lev_r(z;W)-r(z)| W(h)(z)|^2\geq 0\end{equation}
 for $z\in V\setminus \D$ and $W\in \C^{n}$ such that  $W(\rho)(z)=0$.

\noindent \textit{Claim:} To prove (\ref{eqf}) it is  sufficient to prove that for any 
$p\in \D_{\infty}$ we have
 \begin{equation}\label{eqn5}
 \displaystyle \frac{d f(z,W(z))}{d\vec{n}(p)}\bigg|_{z=p}>0  \end{equation} for any smooth vector
field $W$ of type $(1,0)$ such that $W(\rho)=0$ on a neighborhood of $\D_{\infty}$ and $(p,W(p)) \in
 \gd .$ 
 
\noindent \textit{Proof of Claim:} Let us fix $q\in \D_{\infty}$.
Using translation and rotation we can move $q$ to the origin such that the $y_n$-axis is the outward
normal direction at $0.$ There exists a neighborhood $\widetilde{U}$ of $0$ on which $\frac{\partial
\rho}{\partial z_n}$ does not vanish. If $W=(w_1,\ldots ,w_n)$ is a complex tangential vector to the
level set of $\rho$  at $z\in \widetilde{U}$ (i.e. $W(\rho)(z)=0$)  then 
\begin{equation}\label{eqnw_n}w_n=-\left( \left(\frac{\partial \rho(z)}{\partial
z_n}\right)^{-1}\right) \sum_{j=1}^{n-1}\frac{\partial \rho(z)}{\partial z_j}w_j
.\end{equation}

We introduce an auxiliary real valued function $g$  as  
$g(z,W^{\prime} )= f(z,W),$
where $W=(w_1,\ldots ,w_n)$ and $ W^{\prime} =(w_1,\ldots ,w_{n-1}),$ with $w_n$
given by  (\ref{eqnw_n}). 
We choose an open
neighborhood $U$ of $0$ such that $U\Subset  \widetilde{U}.$ Let
$S=\{W^{\prime}\in\C^{n-1}: \|W^{\prime}\|=1\} .$  Notice that $(\overline{U\cap
b\D})\times S$ is compact and it is enough to show that for every $(p,W_{p}^{\prime})\in
(\overline{U\cap b\D})\times S$
there exists a neighborhood $U_p$ of $(p,W_{p}^{\prime})$ in $(\widetilde{U}\cap b\D)\times S$ such
that $g(z,W^{\prime})\geq 0$ for $(z,W^{\prime})\in (U_p\setminus \D)\times S$ and $W^{\prime}\in
\C^{n-1}.$ 
Due to the continuity of the complex Hessian this is true for strongly
pseudoconvex directions. However, (\ref{eqn5}) implies that this is also true for   weakly
pseudoconvex directions. So the proof of the claim is complete.

Let us differentiate $f(z,W(z))$ with respect to $r(z)$ at $p\in \D_{\infty}$. Using (\ref{eqf})  we
get: 
\begin{equation}\label{eqbazuki}  \Lev_h(p;W)+C_r(p;W)+2 Re\bigg( \Lev_r(p;W,d\overline W/dr)\bigg)
- | W(h)(p)|^2 .\end{equation} 
 Since $W$ is a weakly pseudoconvex direction we only need to compute the
complex normal component of $ \frac{d W}{dr(p)}$ at $p$ to estimate the third
term of the above expression. Hence, we need to compute $\frac{d
W}{dr(p)}(r)(p)$
which represents the following: first we differentiate $W$ by $ \frac{d}{dr(p)}$
at $p$ then apply the result to $r$ and evaluate at $p$. Now we use the same calculations used to
derive \eqref{forC_r} and differentiate
the left hand side of $ W(r)(z)+ r(z)W(h)(z)=0$ to get
\begin{eqnarray*}
\frac{d \{ W(r)(z)+ r(z)W(h)(z)\} }{dr(p)}\bigg|_{z=p}=\frac{d W}{dr(p)}(r)(p)+W( \ln \| \nabla r
\|)(p)+ W(h)(p).
\end{eqnarray*} 
 Thus we have:  
 \begin{equation}\label{grad=1}\frac{d W}{dr(p)}(r)(p)+W( \ln \| \nabla r
\|)(p)+ W(h)(p)=0.\end{equation}  
If $Y=\tau N_r+\xi $ where $\xi$ is the complex tangential component of $Y$ then $$\displaystyle
\tau=\frac{Y(r)(p)}{N_r(r)(p)}= Y(r)(p).$$ 
Then using the above observation with (\ref{grad=1}) we conclude that the third term in
(\ref{eqbazuki}) is equal to 
$$-2{\rm Re}\bigg(D_r(p;W)\bigg(\overline{W}(\ln \| \nabla r \|)(p)+ 
\overline{W}(h)(p)\bigg)\bigg).$$
Hence by \eqref{E_rho} we have  $ E_{\rho}(p;W)>0$ for $(p,W(p))\in \gd .$  

\noindent \textit{Modification Away From Infinite Type Points.} In the first part of the proof we
showed that there exist a defining function $\rho$ for $\D$  such that $E_{\rho}(z;W)>0$ for
$(z,W)\in \gd$. That is, there is a neighborhood $V$ of $\D_{\infty}$ such that 
the level sets of $\rho$ are strongly pseudoconvex (from $\D$'s side) in $V\setminus \Dc$.  Now we
will modify $\rho$ away from infinite type points to get a smooth defining function $r$ that will
satisfy i). Let $\rho_{\lambda}(z)=e^{\lambda\rho(z)}-1.$ One can show that 
\begin{equation}\label{eqGlue1}
E_{\rho_{\lambda}}(z;W)=\lambda\Lev_{\rho}(z;W)+E_{\rho}(z;W).
\end{equation}
Then we can choose open sets $V_{1},V_{2}, V_{3},$  and  $\lambda>1$ so that  $\D_{\infty}\Subset
V_{1}\Subset V_{2}\Subset V_{3}\Subset V$ and $E_{\rho_{\lambda}}(z;W)>0$ for $z\in V\cap b\D$ and
$W(\rho)(z)=0 .$ Let $\chi$ be a smooth increasing convex function on the real line so that 
$\chi(t)\equiv 0$ for $t\leq 0$ and $\chi(t) >0$ for $t>0.$ Let us choose 
$A=\sup\{2+\chi '(t):0\leq t\leq 2\} .$ 
Using Lemma \ref{LocalP} we can  choose a sequence of functions $\phi_{j}\in C^{\infty}(\Dc)$ such
that: 
\begin{itemize} 
\item[1)]$-6\ln A< \phi_{j} <- \ln A$ on $\Dc,$
\item[2)] $-6\ln A <\phi_j<-5\ln A$ on $V_{1}\cap b\D,$
\item[3)] $-6\ln A<\phi_j<-3\ln A$ on $V_2 \cap b\D, $
\item [4)] $-2\ln A< \phi_j<-\ln A$ on $b\D\setminus V_3,$ and 
\item[5)] $\Lev_{\phi_{j}}(z;W)> jA^{6}\|W\|^{2}$ for $z\in b\D\setminus V_{1}$ and $W\in \C^{n}$. 
\end{itemize}
Let $h_{j}=e^{\phi_{j}}-1/2$. Then one  can check that $\Lev_{h_j}(z;W)> A|W(h_j)(z)|^2+j$ for $z\in
b\D\setminus V_1$ and $W\in \C^n$. Let us choose
$a=\frac{1}{A^{3}}-\frac{1}{2},\chi_{a}(t)=\chi(t-a),$ and $\psi_j(z)=\chi_{a}\circ h_j(z).$ Then
$\psi_j\equiv 0$ in a neighborhood of $\overline{V_{2}\cap b\D}$. One can calculate that 
\begin{eqnarray*}
 \Lev_{\psi_j}(z;W)&=&\chi_{a}'(h_j(z))\Lev_{h_j}(z;W)+ \chi_{a}''(h_j(z))|W(h_j)(z)|^2 ,\\
|W(\psi_j)(z)|^2&=&|\chi_{a}'(h_j(z))|^2|W(h_j)(z)|^2 .
\end{eqnarray*}
Let $r(z)=\rho_{\lambda}(z) e^{\psi_j(z)}$. As in \eqref{E_rho} one can calculate that 
$$E_{r}(z;W)=\Lev_{\psi_j}(z;W)-|W(\psi_j)(z)|^2-2{\rm
Re}\Big(D_{\rho_{\lambda}}(z;W)\overline{W}(\psi_j)(z)\Big)+E_{\rho_{\lambda}}(z;W) .$$
Since  $r(z)=\rho_{\lambda}(z)$ in a neighborhood of $\overline{V_{2}\cap b\D}$ we only need to show
that 
\begin{eqnarray}\label{eqglue}
\Lev_{h_j}(z;W)
>(\chi_{a}'(h_j(z))+1)|W(h_j(z))|^2+|D_{\rho_{\lambda}}(z;W)|^2-\frac{E_{\rho_{\lambda}}(z;W)}{\chi_
{a}'(h_j)}
\end{eqnarray}
for $z\in b\D \setminus V_{2}$ and $W(r)(z)=0 .$
Let us choose 
\begin{equation}\label{eqGlue}j>\sup\left\{|D_{\rho_{\lambda}}(z;W)|^2-\frac{E_{\rho_{\lambda}}(z;W)
}{\chi_{a}'(h_k(z))}:z\in b\D, W(r)(z)=0,k=1,2,\ldots \right\} .\end{equation}
We note that $\chi_{a}$ and $\phi_{j}$'s are chosen so that $\chi_{a}'\geq 0$ and
$E_{\rho_{\lambda}}(z;W)>0$ for $z\in V\cap b\D$ and $W(\rho)(z)=0 $.  $E_{\rho_{\lambda}}$ can be
negative outside $V$ in some directions but  there exists $b>0$ such that $\chi_{a}'(h_j(z))>b$ for
$j=1,2,\ldots$ and $z\in b\D \setminus V_{3}$. So the right hand side of  (\ref{eqGlue}) is finite
and, since $A=\sup\{2+\chi '(t):0\leq t\leq 2\}$,  one can choose $j$ so that \eqref{eqglue} is
satisfied. Hence we showed that 
$E_{r}(z;W)> 0$ for $z\in b\D$ and $W(r)(z)=0$. Similar argument used in the proof of
ii)$\Rightarrow$iii) shows that  there exists a neighborhood $U$ of $\Dc$ and $c>0$ such that
$\Lev_{r}(z;W)\geq cr(z)\|W\|^{2}$ for $z\in U\setminus \D$ and $W(r)(z)=0.$
\end{proof}

Now we will give the proof of Corollary \ref{Invariance}.

\begin{proof}[Proof of Corollary \ref{Invariance}] 
Let $n_{1}(p)$ denote the unit outward normal of $\D_{1}$ at a boundary point $p\in b\D_{1}$.
Furthermore, let $q=F(p)$ and $n_{2}=F(n_{1})$. We note that $n_{2}$ is transversal to $b\D_{2}$ at
$q$ and i) in Theorem \ref{MainThm} implies that there exist a defining function $\rho_{2}$ for 
$\D_{2}$  and $c_{2}>0$ such that
$$\frac{d}{dn_{2}(q)}\Lev_{\rho_{2}}(z;W(z))\bigg|_{z=q} >c_{2}\|W\|^{2} $$  
for $q\in b\D_{2}$ and $W(\rho_{2})=0$ in a neighborhood of $\Dc_{2}$. Namely,
\begin{equation}\label{ineq1}
\lim_{\varepsilon\to 0^{+}}\frac{\Lev_{\rho_{2}}(q-\varepsilon n_{2}(q);W(q-\varepsilon
n_{2}(q)))}{\varepsilon}>c_{2} 
\end{equation}
for $(q,W(q))\in \gdtwo $.
Let $\rho_{1}(z)=\rho_{2}(F(z))$ and  extend $F$ smoothly to some neighborhood of $\Dc_{1}$. So
$\rho_{1}$ is a defining function for $\D_{1}.$ Since $F$ is proper it transforms the complex
Hessian of $\rho_{2}$ to the complex Hessian of $\rho_{1}$. More precisely, let $J_{F}$ denote the
complex Jacobian of $F$. That is, 
$$J_{F}=\left\{\frac{\partial F_{j}}{\partial z_{k}} \right\}_{j,k} .$$ Then one can show that
$\Lev_{\rho_{1}}(z;W(z))=\Lev_{\rho_{2}}(F(z),J_{F}W(F(z)))$ for $z\in \Dc_{1}$ and
$W(\rho_{1})(z)=0$ if an only if $J_{F}W(\rho_{2})(F(z))=0.$ Let us fix a smooth vector field
$\widetilde{W}$ of type $(1,0)$ that is complex tangential to level sets of $\rho_{1}$ on $\Dc_{1}$
and denote $W=J_{F}\widetilde{W}$.  Since $F$ is holomorphic and extends to the boundary
(\ref{ineq1}) implies that  
$$ \lim_{\varepsilon\to 0^{+}}\frac{\Lev_{\rho_{1}}(p-\varepsilon
n_{1}(p);\widetilde{W}(p-\varepsilon n_{1}(p)))}{\varepsilon}\geq 0  $$
for $(p,\widetilde{W}(p))\in \gdone $.
So we have $E_{\rho_{1}}(z;A)\geq 0$ for $(z,A)\in \gdone .$ Therefore, Theorem \ref{MainThm}
implies that $\Dc_{1}$ has a strong Stein neighborhood basis. 
\end{proof}

\section{Proof of Corollary \ref{P-tilde}}\label{sec4}

\begin{lemma}\label{LocalGlobal} Let $\D$ be a smooth bounded  pseudoconvex domain in $\C^n.$ Assume
that for every $M>0$ there exist a neighborhood $U$ of 
 $b\D$ and $\psi\in C^{\infty}(U)$ such that $ |W(\psi)(z)|^2\leq \Lev_{\psi}(z;W)$ and
$\Lev_{\psi}(z;W)\geq M\|W\|^{2},$ for 
 $z\in U$ and $W\in \C^{n}.$ Then $\D$ satisfies property $\ptilde$.
\end{lemma}
\begin{proof} 

One can check that $\Lev_h(z;W)\geq|W(h)(z)|^2$ if and only if $-e^{-h}$ is plurisubharmonic at $z$
in the direction $W$. Without loss of generality we may assume that $\psi\leq -1 $ on $U$. Let
$A=-\sup_{z\in b\D}e^{-\psi(z)}$. Since $\D\setminus U$ is compact we can use a theorem of Demailly
\cite{De87} to choose Green's functions $G_1,\cdots,G_k$ so that $G_1+\cdots+G_k-e^{-\psi}<A-1$ on
$\D\setminus U$. Then  $F(z)=\max\{G_1(z)+\cdots+G_k(z)-e^{-\psi(z)}, A-1/2\}$ is a continuous
plurisubharmonic function that can be extended to a neighborhood of $\Dc$. Using convolution with an
approximate identity we may choose a function $\widetilde F\in C^{\infty}(\Dc)$ that is
plurisubharmonic on $\D$ and arbitrarily close to $F$ uniformly on $\Dc$. Let $f(z)=-\log(-F(z))$ 
and $\tilde f(z)=-\log( -\widetilde F(z))$. Since $\widetilde F$ is plurisubharmonic we have
$\Lev_{\tilde f}(z;W)\geq|W(\tilde f)(z)|^2$ for $z\in \Dc$ and $W\in \C^n$. Now we need to show
that the complex Hessian of $\tilde f$ is large enough on $b\D$. Since on a sufficiently small band
close to $b\D$ we have $F(z)=G_1(z)+\cdots+G_k(z)-e^{-\psi(z)}$ one can show that
\begin{eqnarray*}
 \Lev_{f}(z;W)&=&\frac{\Lev_{g}(z;W)+e^{-\psi(z)}(\Lev_{\psi}(z;W)-|W(\psi)(z)|^{2})}{e^{-\psi(z)}
-g(z)}\\ &&+\frac{|W(g)(z)+e^{-\psi(z)}W(\psi)(z)|^{2}}{(e^{-\psi(z)}-g(z))^{2}}
\end{eqnarray*}
 where $g(z)=G_{1}(z)+\cdots+G_{k}(z)$. Since $g\equiv 0$ on $b\D$ we have $\Lev_{f}(z;W)\geq
\Lev_{\psi}(z;W)-|W(\psi)(z)|^{2}$ on $b\D$. But we could have chosen $\psi$ so that
$\Lev_{\psi}(z;W)\geq 2|W(\psi)(z)|^{2}$ which would imply that $\Lev_{f}(z;W)\geq
\frac{1}{2}\Lev_{\psi}(z;W)$ on $b\D$. We can choose $\tilde f$ sufficiently close to $f$ so that
$\Lev_{\tilde f}(z;W)\geq \frac{1}{2}\Lev_{\psi}(z;W)-\|W\|^2$ for $z\in b\D$ and $W\in \C^n$.
Hence, $\D$ satisfies property $\ptilde$.
\end{proof} 
\begin{proof}[Proof of Corollary \ref{P-tilde}] Using the proof of Theorem \ref{MainThm} one can
easily prove that if $\D$ satisfies property $\ptilde$ then for every $M>0$ there exists a defining
function $\rho$ such that $E_{\rho}(z;W)>M$ for $(z;W)\in \gd .$  To prove the other direction let
us assume that for $M>0$ there exists a defining function $\rho$ such that $E_{\rho}(z;W)>M$ for
$(z,W)\in \gd$. Let us define $\rho_{\lambda}(z)=e^{\lambda \rho(z)}-1$ and fix a defining function
$r$ for $\D$. Then by (\ref{eqGlue1}) we can choose a large enough $\lambda>1$ so that
$E_{\rho_{\lambda}}(z;W)>M\|W\|^{2}$ for $z\in b\D$ and $W(\rho)(z)=0.$ Then there exist $h\in
C^{\infty}(\Dc)$ such that  $\rho_{\lambda}(z)=r(z)e^{h(z)}$. By \eqref{E_rho} the condition
$E_{\rho_{\lambda}}(z;W)>M$ implies that 
\begin{eqnarray*}
\Lev_h(z;W)&>&M\|W\|^2+|W(h)(z)|^2+2{\rm Re}\bigg(D_r(z;W)\overline{W}(h)(z)\bigg)-E_{r}(z;W)\\
&>&M\|W\|^2+\frac{|W(h)(z)|^2}{2}-2|D_r(z;W)|^{2}-|E_{r}(z;W)|
\end{eqnarray*}
for $z\in b\D$ and $W\in\C^n .$ Let $\tilde h(z)=h(z)/2$ and  
$$\widetilde M=\frac{M}{2}-\sup\bigg\{|D_r(z;W)|^{2}+\frac{1}{2}|E_{r}(z;W)|:(z,W)\in b\D\times\C^n,
\|W\|\leq 1\bigg\}.$$ Since we can choose $M$ as large as we wish for  every $\widetilde M>0$ there
exist a neighborhood $U$ of $b\D$ and  $\tilde h\in C^{\infty}(U)$  such that $\Lev_{\tilde
h}(z;W)>\widetilde M+|W(\tilde h)|^{2}$ for $z\in U$ and $W\in \C^{n}$.  Then Lemma
\ref{LocalGlobal} implies that 
$\D$ satisfies property $\ptilde$.
\end{proof}
     
\singlespace


\begin{thebibliography}{McN02}

\bibitem[BF78]{BF78}
Eric Bedford and John~Erik Forn{\ae}ss, \emph{Domains with pseudoconvex
 neighborhood systems}, Invent. Math. \textbf{47} (1978), no.~1, 1--27,
available at {\tt http://gdz.sub.uni-goettingen.de}.

\bibitem[Bell86]{Bell86}
Steve~Bell, \emph{Differentiability of the {B}ergman kernel and pseudolocal estimates},
Math. Z. \textbf{192} (1986), no,~3 467--472, 
available at {\tt http://gdz.sub.uni-goettingen.de}.

\bibitem[BS91]{BS91}
Harold~P. Boas and Emil~J. Straube, \emph{Sobolev estimates for the
 {$\overline\partial$}-{N}eumann operator on domains in {${\bf C}\sp n$}
 admitting a defining function that is plurisubharmonic on the boundary},
 Math. Z. \textbf{206} (1991), no.~1, 81--88, available at {\tt
 http://gdz.sub.uni-goettingen.de}.

\bibitem[BS93]{BS93}
\bysame, \emph{de {R}ham cohomology of manifolds containing the points of
 infinite type, and {S}obolev estimates for the
{$\overline\partial$}-{N}eumann problem}, J. Geom. Anal. \textbf{3} (1993),
 no.~3, 225--235.

\bibitem[BS99]{BS}
\bysame, \emph{Global regularity of the
 {$\overline\partial$}-{N}eumann problem: a survey of the {$L\sp 2$}-{S}obolev
 theory}, Several complex variables (Berkeley, CA, 1995--1996), Math. Sci.
 Res. Inst. Publ., vol.~37, Cambridge Univ. Press, Cambridge, 1999,
  pp.~79--111, arXiv:math.CV/9612204.

\bibitem[BT33]{BehThullen}
H.~Behnke and P.~Thullen, \emph{Zur {T}heorie der {S}ingularit\"{a}ten der
 {F}unktionen mehrerer komplexen {V}er\"{a}nderlichen. {D}as
 {K}onvergenzproblem der {R}egularit\"{a}tsh\"{u}llen}, Math. Ann.
 \textbf{108} (1933), 91--104, available at {\tt
 http://gdz.sub.uni-goettingen.de}.

\bibitem[Cat84]{Cat}
David~W. Catlin, \emph{Global regularity of the {$\bar \partial $}-{N}eumann
  problem}, Complex analysis of several variables (Madison, Wis., 1982), Proc.
  Sympos. Pure Math., vol.~41, Amer. Math. Soc., Providence, RI, 1984,
  pp.~39--49.
  
\bibitem[\c{C}e]{Celik}
Mehmet \c{C}elik, \emph{Contributions to the compactness theory of the  {$\bar
\partial $}-{N}eumann operator,} Ph.D. thesis, Texas A{\&}M University, TX,
2008,  available at {\tt http://hdl.handle.net/1969.1/ETD-TAMU-2008-05-6}

\bibitem[{\v{C}}ir69]{Cir}
E.~M. {\v{C}}irka, \emph{Approximation by holomorphic functions on smooth
 manifolds in {${\bf C}\sp{n}$}}, Mat. Sb. (N.S.) \textbf{78 (120)} (1969),
 101--123, Math. USSR Sb., \textbf{7} (1969), 95--113.

\bibitem[D82]{D} John~P.~ D'Angelo, \textit{Real hypersurfaces, orders of contact, and
applications}, Ann. of Math. (2) \textbf{115} (1982), no.~3, 615--637.

\bibitem[De87]{De87} Jean-Pierre~ Demailly, \textit{Mesures de {M}onge-{A}mp\`ere et mesures
pluriharmoniques}, Math. Z.  \textbf{194}(1987),  no.~ 4, 519--564.
   
\bibitem[DF77a]{DF77a}
Klas Diederich and John~Erik Forn{\ae}ss, \emph{Pseudoconvex domains: an
  example with nontrivial {N}ebenh\"ulle}, Math. Ann. \textbf{225} (1977),
  no.~3, 275--292, available at {\tt http://gdz.sub.uni-goettingen.de}.

\bibitem[DF77b]{DF77b}
\bysame, \emph{Pseudoconvex domains: existence of {S}tein neighborhoods}, Duke
 Math. J. \textbf{44} (1977), no.~3, 641--662.

\bibitem[FH07]{FH}
John~Erik Fornaess and Anne-Katrin Herbig, \emph{{A note on plurisubharmonic
 defining functions in $\mathbb{C}^2$}},  Math. Z.  \textbf{257} (2007),  no.
4, 769--781.

\bibitem[FH08]{FH2}
John~Erik Fornaess and Anne-Katrin Herbig, \emph{{A note on plurisubharmonic
 defining functions in $\mathbb{C}^n$}},   Math. Ann.  \textbf{342} (2008), 
no. 4, 749--772.

\bibitem[FN77]{FN}
John~Erik Forn{\ae}ss and Alexander Nagel, \emph{The {M}ergelyan property for
 weakly pseudoconvex domains}, Manuscripta Math. \textbf{22} (1977), no.~2,
 199--208, available at {\tt http://gdz.sub.uni-goettingen.de}.

\bibitem[Har06]{Har06}
Phillip~S. Harrington, \emph{Property $(P)$ and Stein neighborhood 
bases on $C^1$ domains},  Illinois J. Math.  \textbf{52} (2008),  no. 1,
145--151


\bibitem[H{\"o}r90]{Ho}
Lars H{\"o}rmander, \emph{An introduction to complex analysis in several
  variables}, third ed., North-Holland Mathematical Library, vol.~7,
  North-Holland Publishing Co., Amsterdam, 1990.
  
\bibitem[Kra01]{Kr}
Steven~G. Krantz, \emph{Function theory of several complex variables}, AMS
  Chelsea Publishing, Providence, RI, 2001, Reprint of the 1992 edition.
  
\bibitem[McN02]{McN}
Jeffery~D. McNeal, \emph{A sufficient condition for compactness of the
 {$\overline\partial$}-{N}eumann operator}, J. Funct. Anal. \textbf{195}
 (2002), no.~1, 190--205.

\bibitem[Ran86]{Ra}
R.~Michael Range, \emph{Holomorphic functions and integral representations in
  several complex variables}, Graduate Texts in Mathematics, vol. 108,
  Springer-Verlag, New York, 1986.


\bibitem[Ros61]{Ros61}
Hugo Rossi, \emph{Holomorphically convex sets in several complex variables},
 Ann. of Math. (2) \textbf{74} (1961), 470--493.

\bibitem[Sib87a]{Sib87}
Nessim Sibony, \emph{Une classe de domaines pseudoconvexes}, Duke Math. J.
 \textbf{55} (1987), no.~2, 299--319.

\bibitem[Sib87b]{Sib87a} \bysame,\emph{Hypoellipticit\'e pour l'op\'erateur {$\overline\partial$}},
Math. Ann. \textbf{276}  (1987), no.~2, 279--290, available at {\tt
http://gdz.sub.uni-goettingen.de}. 

\bibitem[Sib91]{Sib91}
\bysame, \emph{Some aspects of weakly pseudoconvex domains}, Several complex
 variables and complex geometry, Part 1 (Santa Cruz, CA, 1989), Proc. Sympos.
Pure Math., vol.~52, Amer. Math. Soc., Providence, RI, 1991, pp.~199--231.


\bibitem[Ste87]{Sten}
Berit Stens{\o}nes, \emph{Stein neighborhoods}, Math. Z. \textbf{195} (1987),
 no.~3, 433--436, available at {\tt http://gdz.sub.uni-goettingen.de}.

\bibitem[Str01]{St01}
Emil~J. Straube, \emph{Good {S}tein neighborhood bases and regularity of the
 {$\overline\partial$}-{N}eumann problem}, Illinois J. Math. \textbf{45}
(2001), no.~3, 865--871, arXiv: math.CV/0006136.

\bibitem[\c{S}ah06]{Sah06}
S\"{o}nmez \c{S}ahuto\u{g}lu, \emph{Compactness of the
$\overline{\partial}$-{N}eumann problem and {S}tein neighborhood bases}, Ph.D.
thesis, Texas A{\&}M University, TX, 2006, available at 
{\tt http://handle.tamu.edu/1969.1/3879}.

\end{thebibliography}
\end{document}